\documentclass[a4paper,12pt]{article}
\usepackage{latexsym}
\usepackage{amsfonts}
\usepackage{eucal}
\usepackage{amsmath}

\numberwithin{equation}{section}

\input{tcilatex}

\begin{document}

\title{\textbf{Subclasses of meromorphically multivalent functions defined
by a differential operator}}
\date{}
\author{\textbf{\small Halit Orhan}, \textbf{\small Dorina R\u aducanu} 
{\small and \textbf{Erhan Deniz}}}
\maketitle

\abstract
In this paper we introduce and study two new subclasses $\Sigma_{\lambda\mu
mp}(\alpha,\beta)$ and $\Sigma^{+}_{\lambda\mu mp}(\alpha,\beta)$ of
meromorphically multivalent functions which are defined by means of a new
differential operator. Some results connected to subordination properties,
coefficient estimates, convolution properties, integral representation,
distortion theorems are obtained. We also extend the familiar concept of $%
(n,\delta)-$neighborhoods of analytic functions to these subclasses of
meromorphically multivalent functions. \endabstract

\textit{2000 Mathematics Subject Classification:} 30C45.

\textit{Key words and phrases:} Analytic functions, meromorphic functions,
multivalent functions, differential operator, subordination, neighborhoods.

\newtheorem{theo}{Theorem}[section] \newtheorem{cor}{Corollary}[section] %
\newtheorem{defn}{Definition}[section] \newtheorem{rem}{Remark}[section] %
\newtheorem{lem}{Lemma}[section]

\section{Introduction}

Let $\tilde{\mathcal{A}}$ be the class of analytic functions in the unit
disk $\mathbb{U}=\left\{z\in\mathbb{C}:|z|<1\right\}$.

Consider 
\begin{equation}  \label{1.1}
\Omega=\left\{w\in\tilde{\mathcal{A}}:w(0)=0\;\text{and}\;|w(z)|<1, z\in%
\mathbb{U}\right\}
\end{equation}
the class of Schwarz functions.

For $0\leq\alpha<1$ let 
\begin{equation}  \label{1.2}
\mathcal{P}(\alpha)=\left\{p\in\tilde{\mathcal{A}}:p(0)=1\;\text{and}\;\Re
p(z)>\alpha, z\in\mathbb{U}\right\}.
\end{equation}
Note that $\mathcal{P}=\mathcal{P}(0)$ is the well-known Carath$\acute{e}$%
odory class of functions.

The classes of Schwarz and Carath$\acute{e}$odory functions play an
extremely important role in the theory of analytic functions and have been
studied by many authors.

It is easy to see that 
\begin{equation}  \label{1.3}
p\in\mathcal{P}(\alpha)\;\;\text{if and only if}\;\;\frac{p(z)-\alpha}{%
1-\alpha}\in\mathcal{P}.
\end{equation}
Making use of the properties of functions in the class $\mathcal{P}$ and the
condition (1.3), the following properties of the functions in the class $%
\mathcal{P}(\alpha)$ can be obtained. 
\begin{lem}
Let $p\in\tilde{\mathcal{A}}$. Then $p\in\mathcal{P}(\alpha)$ if and only if there exists $w\in\Omega$ such that
\begin{equation}\label{1.4}
p(z)=\frac{1-(2\alpha-1)w(z)}{1-w(z)}\;\;(z\in\mathbb{U}).
\end{equation}
\end{lem}
\begin{lem}(Herglotz formula)
A function $p\in\tilde{\mathcal{A}}$ belongs to the class $\mathcal{P}(\alpha)$ if and only if there exists a probability measure $\mu(x)$ on $\partial\mathbb{U}$ such that
\begin{equation}\label{1.5}
p(z)=\int_{|x|=1}\frac{1-(2\alpha-1)xz}{1-xz}d\mu(x)\;\;(z\in\mathbb{U}).
\end{equation}
The correspondence between $\mathcal{P}(\alpha)$ and probability measure $\mu(x)$ on $\partial\mathbb{U}$, given by (1.5) is one-to-one.
\end{lem}

If $f$ and $g$ are in $\tilde{\mathcal{A}}$, we say that $f$ is subordinate
to $g$, written $f\prec g$, if there exists a function $w\in\Omega$ such
that $f(z)=g(w(z))\;(z\in\mathbb{U})$. It is known that if $f\prec g$, then $%
f(0)=g(0)$ and $f(\mathbb{U})\subset g(\mathbb{U})$. In particular, if $g$
is univalent in $\mathbb{U}$ we have the following equivalence: 
\begin{equation*}
f(z)\prec g(z)\;\;(z\in\mathbb{U})\;\;\text{if and only if}\;\;f(0)=g(0)\;\;%
\text{and}\;\;f(\mathbb{U})\subset g(\mathbb{U}).
\end{equation*}

Let $\Sigma_{p}$ denote the class of all meromorphic functions $f$ of the
form 
\begin{equation}  \label{1.6}
f(z)=z^{-p}+\sum_{k=1-p}^{\infty}a_{k}z^{k}\;\;(p\in\mathbb{N}%
:=\left\{1,2,\ldots\right\})
\end{equation}
which are analytic and $p$-valent in the punctured unit disk\newline
$\mathbb{U}^{*}=\mathbb{U}\setminus\left\{0\right\}$.

Denote by $\Sigma_{p}^{+}$ the subclass of $\Sigma_{p}$ consisting of
functions of the form 
\begin{equation}  \label{1.7}
f(z)=z^{-p}+\sum_{k=1-p}^{\infty}a_{k}z^{k}\;\;a_{k}\geq0\;\;(z\in\mathbb{U}%
^{*}).
\end{equation}

A function $f\in\Sigma_{p}$ is meromorphically multivalent starlike of order 
$\alpha\;(0\leq\alpha<1)$ (see \cite{aou}) if 
\begin{equation*}
-\Re\left\{\frac{1}{p}\frac{zf^{\prime }(z)}{f(z)}\right\}>\alpha\;\;(z\in%
\mathbb{U}).
\end{equation*}
The class of all such functions is denoted by $\Sigma_{p}^{*}(\alpha)$.

If $f\in\Sigma_{p}$ is given by (1.6) and $g\in\Sigma_{p}$ is given by 
\begin{equation*}
g(z)=z^{-p}+\sum_{k=1-p}^{\infty}b_{k}z^{k}
\end{equation*}
then the Hadamard product (or convolution) of $f$ and $g$ is defined by 
\begin{equation*}
(f\ast g)(z)=z^{-p}+\sum_{k=1-p}^{\infty}a_{k}b_{k}z^{k}=(g\ast
f)(z)\;\;(p\in\mathbb{N}\;\;z\in\mathbb{U}^{*}).
\end{equation*}

For a function $f\in\Sigma_{p}$, we define the differential operator $%
D_{\lambda\mu p}^{m}$ in the following way: 
\begin{equation*}
D_{\lambda\mu p}^{0}f(z)=f(z)
\end{equation*}
\begin{equation}  \label{1.8}
D_{\lambda\mu p}^{1}f(z)=D_{\lambda\mu p}f(z)=\lambda\mu\frac{%
[z^{p+1}f(z)]^{\prime \prime }}{z^{p-1}}+(\lambda-\mu)\frac{%
[z^{p+1}f(z)]^{\prime }}{z^{p}}+(1-\lambda+\mu)f(z).
\end{equation}
and, in general 
\begin{equation}  \label{1.9}
D_{\lambda\mu p}^{m}f(z)=D_{\lambda\mu p}(D_{\lambda\mu p}^{m-1}f(z)),
\end{equation}
where $0\leq\mu\leq\lambda$ and $m\in\mathbb{N}$.

If the function $f\in\Sigma_{p}$ is given by (1.6) then, from (1.8) and
(1.9), we obtain 
\begin{equation}  \label{1.10}
D_{\lambda\mu
p}^{m}f(z)=z^{-p}+\sum_{k=1-p}^{\infty}\Phi_{k}(\lambda,\mu,m,p)a_{k}z^{k}
\end{equation}
\begin{equation*}
(m\in\mathbb{N}\;,\;p\in\mathbb{N}\;,\;z\in\mathbb{U}^{*})
\end{equation*}
where 
\begin{equation}  \label{1.11}
\Phi_{k}(\lambda,\mu,m,p)=[1+(k+p)(\lambda-\mu+(k+p+1)\lambda\mu)]^{m}.
\end{equation}

From (1.10) it follows that $D_{\lambda\mu p}^{m}f(z)$ can be written in
terms of convolution as 
\begin{equation}  \label{1.12}
D_{\lambda\mu p}^{m}f(z)=(f\ast h)(z)
\end{equation}
where 
\begin{equation}  \label{1.13}
h(z)=z^{-p}+\sum_{k=1-p}^{\infty}\Phi_{k}(\lambda,\mu,m,p)z^{k}.
\end{equation}

Note that, the case $\lambda=1$ and $\mu=0$ of the differential operator $%
D_{\lambda\mu p}^{m}$ was introduced by Srivastava and Patel \cite{sri}. A
special case of $D_{\lambda\mu p}^{m}$ for $p=1$ was considered in \cite{rad}%
.

Making use of the differential operator $D_{\lambda\mu p}^{m}$, we define a
new subclass of the function class $\Sigma_{p}$ as follows. 
\begin{defn}
A function $f\in\Sigma_{p}$ is said to be in the class $\Sigma_{\lambda\mu m p}(\alpha,\beta)$ if it satisfies the condition
\begin{equation}\label{1.14}
\left|\frac{1}{p}\frac{z(D_{\lambda\mu p}^{m}f(z))'}{D_{\lambda\mu p}^{m}f(z)}+1\right|<\beta\left|\frac{1}{p}\frac{z(D_{\lambda\mu p}^{m}f(z))'}{D_{\lambda\mu p}^{m}f(z)}+2\alpha-1\right|
\end{equation}
for some $\alpha\;(0\leq\alpha<1)$, $\beta\;(0<\beta\leq1)$ and $z\in\mathbb{U}^{*}$.
\end{defn}
Note that a special case of the class $\Sigma_{\lambda\mu m p}(\alpha,\beta)$
for $p=1$ and $m=0$ is the class of meromorphically starlike functions of
order $\alpha$ and type $\beta$ introduced earlier by Mogra et al. \cite{mog}%
. It is easy to check that for $m=0$ and $\beta=1$, the class $%
\Sigma_{\lambda\mu m p}(\alpha,\beta)$ reduces to the class $%
\Sigma_{p}^{*}(\alpha)$.

We consider another subclass of $\Sigma_{p}$ given by 
\begin{equation}  \label{1.15}
\Sigma_{\lambda\mu m
p}^{+}(\alpha,\beta):=\Sigma_{p}^{+}\cap\Sigma_{\lambda\mu m
p}(\alpha,\beta).
\end{equation}

The main object of this paper is to present a systematic investigation of
the classes $\Sigma_{\lambda\mu m p}(\alpha,\beta)$ and $\Sigma_{\lambda\mu
m p}^{+}(\alpha,\beta)$.

\section{Properties of the class $\Sigma_{\protect\lambda\protect\mu m p}(%
\protect\alpha,\protect\beta)$}

We begin this section with a necessary and sufficient condition, in terms of
subordination, for a function to be in the class $\Sigma_{\lambda\mu m
p}(\alpha,\beta)$. 
\begin{theo}
A function $f\in\Sigma_{p}$ is in the class $\Sigma_{\lambda\mu m p}(\alpha,\beta)$ if and only if
\begin{equation}\label{2.1}
\frac{z(D_{\lambda\mu p}^{m}f(z))'}{D_{\lambda\mu p}^{m}f(z)}\prec\frac{p(2\alpha-1)\beta z-p}{1-\beta z}\;(z\in\mathbb{U}).
\end{equation}
\end{theo}
\textit{Proof.} Let $f\in\Sigma_{\lambda\mu m p}(\alpha,\beta)$. Then, from
(1.6), we have 
\begin{equation*}
\left|-\frac{1}{p}\frac{z(D_{\lambda\mu p}^{m}f(z))^{\prime }}{D_{\lambda\mu
p}^{m}f(z)}-1\right|^{2}<\beta^{2}\left|-\frac{1}{p}\frac{z(D_{\lambda\mu
p}^{m}f(z))^{\prime }}{D_{\lambda\mu p}^{m}f(z)}+1-2\alpha\right|^{2}
\end{equation*}
or 
\begin{equation*}
(1-\beta^{2})\left|-\frac{1}{p}\frac{z(D_{\lambda\mu p}^{m}f(z))^{\prime }}{%
D_{\lambda\mu p}^{m}f(z)}\right|^{2}-2[1+\beta^{2}(1-2\alpha)]\Re\left\{-%
\frac{1}{p}\frac{z(D_{\lambda\mu p}^{m}f(z))^{\prime }}{D_{\lambda\mu
p}^{m}f(z)}\right\}
\end{equation*}
\begin{equation*}
<\beta^{2}(1-2\alpha)^{2}-1.
\end{equation*}
If $\beta\neq1$, we have 
\begin{equation*}
\left|-\frac{1}{p}\frac{z(D_{\lambda\mu p}^{m}f(z))^{\prime }}{D_{\lambda\mu
p}^{m}f(z)}\right|^{2}-2\frac{1+\beta^{2}(1-2\alpha)}{1-\beta^{2}}\Re\left\{-%
\frac{1}{p}\frac{z(D_{\lambda\mu p}^{m}f(z))^{\prime }}{D_{\lambda\mu
p}^{m}f(z)}\right\}
\end{equation*}
\begin{equation*}
+\left[\frac{1+\beta^{2}(1-2\alpha)}{1-\beta^{2}}\right]^{2}<\frac{%
\beta^{2}(2\alpha-1)^{2}-1}{1-\beta^{2}}+\left[\frac{1+\beta^{2}(1-2\alpha)}{%
1-\beta^{2}}\right]^{2},
\end{equation*}
that is 
\begin{equation}  \label{2.2}
\left|-\frac{1}{p}\frac{z(D_{\lambda\mu p}^{m}f(z))^{\prime }}{D_{\lambda\mu
p}^{m}f(z)}-\frac{1-\beta^{2}(2\alpha-1)}{1-\beta^{2}}\right|<\frac{%
2\beta(1-\alpha)}{1-\beta^{2}}.
\end{equation}
The inequality (2.2) shows that the values region of $F(z)=-\frac{1}{p}\frac{%
z(D_{\lambda\mu p}^{m}f(z))^{\prime }}{D_{\lambda\mu p}^{m}f(z)}$ is
contained in the disk centered at $\frac{1-\beta^{2}(2\alpha-1)}{1-\beta^{2}}
$ and radius $\frac{2\beta(1-\alpha)}{1-\beta^{2}}$. It is easy to check
that the function $G(z)=\frac{1-(2\alpha-1)\beta z}{1-\beta z}$ maps the
unit disk $\mathbb{U}$ onto the disk 
\begin{equation*}
\left|\omega-\frac{1-\beta^{2}(2\alpha-1)}{1-\beta^{2}}\right|<\frac{%
2\beta(1-\alpha)}{1-\beta^{2}}.
\end{equation*}
Since $G$ is univalent and $F(0)=G(0)$, $F(\mathbb{U})\subset G(\mathbb{U})$%
, we obtain that $F(z)\prec G(z)$, that is 
\begin{equation*}
-\frac{1}{p}\frac{z(D_{\lambda\mu p}^{m}f(z))^{\prime }}{D_{\lambda\mu
p}^{m}f(z)}\prec\frac{1-(2\alpha-1)\beta z}{1-\beta z}
\end{equation*}
or 
\begin{equation*}
\frac{z(D_{\lambda\mu p}^{m}f(z))^{\prime }}{D_{\lambda\mu p}^{m}f(z)}\prec%
\frac{p(2\alpha-1)\beta z-p}{1-\beta z}.
\end{equation*}
Conversely, suppose that subordination (2.1) holds. Then 
\begin{equation}  \label{2.3}
-\frac{1}{p}\frac{z(D_{\lambda\mu p}^{m}f(z))^{\prime }}{D_{\lambda\mu
p}^{m}f(z)}=\frac{1-(2\alpha-1)\beta w(z)}{1-\beta w(z)}
\end{equation}
where $w\in\Omega$. After simple calculations, from (2.3), we obtain 
\begin{equation*}
\left|\frac{1}{p}\frac{z(D_{\lambda\mu p}^{m}f(z))^{\prime }}{D_{\lambda\mu
p}^{m}f(z)}+1\right|<\beta\left|\frac{1}{p}\frac{z(D_{\lambda\mu
p}^{m}f(z))^{\prime }}{D_{\lambda\mu p}^{m}f(z)}+2\alpha-1\right|
\end{equation*}
which proves that $f\in\Sigma_{\lambda\mu m p}(\alpha,\beta)$.

If $\beta=1$, inequality (1.14) becomes 
\begin{equation}  \label{2.4}
\left|-\frac{1}{p}\frac{z(D_{\lambda\mu p}^{m}f(z))^{\prime }}{D_{\lambda\mu
p}^{m}f(z)}-1\right|<\left|-\frac{1}{p}\frac{z(D_{\lambda\mu
p}^{m}f(z))^{\prime }}{D_{\lambda\mu p}^{m}f(z)}+1-2\alpha\right|.
\end{equation}
From (2.4) we can easily obtain that 
\begin{equation*}
-\frac{1}{p}\frac{z(D_{\lambda\mu p}^{m}f(z))^{\prime }}{D_{\lambda\mu
p}^{m}f(z)}\prec \frac{1-(2\alpha-1)z}{1-z}
\end{equation*}
or 
\begin{equation}  \label{2.5}
\frac{z(D_{\lambda\mu p}^{m}f(z))^{\prime }}{D_{\lambda\mu p}^{m}f(z)}\prec 
\frac{(2\alpha-1)pz-p}{1-z}.
\end{equation}
Thus, the proof of our theorem is completed. 
\begin{rem}
Since $\Re\displaystyle\frac{1-(2\alpha-1)\beta z}{1-\beta z}>\alpha$ it follows that
$$-\Re\left\{\frac{1}{p}\frac{z(D_{\lambda\mu p}^{m}f(z))'}{D_{\lambda\mu p}^{m}f(z)}\right\}>\alpha$$
which shows that $D_{\lambda\mu p}^{m}f(z)\in\Sigma_{p}^{*}(\alpha)$.
\end{rem}

Making use of the subordination relationship for the class $%
\Sigma_{\lambda\mu m p}(\alpha,\beta)$, we derive a structural formula,
first for the class $\Sigma_{\lambda\mu m p}(\alpha,1)$ and then for the
class $\Sigma_{\lambda\mu m p}(\alpha,\beta)$. 
\begin{theo}
A function $f\in\Sigma_{p}$ belongs to the class $\Sigma_{\lambda\mu m p}(\alpha,1)$ if and only if there exists a probability measure $\mu(x)$ on $\partial\mathbb{U}$ such that
$$f(z)=\left[z^{-p}+\sum_{k=1-p}^{\infty}\frac{z^{k}}{\Phi_{k}(\lambda,\mu,m,p)}\right]$$
\begin{equation}\label{2.6}
\ast\left[z^{-p}\exp\int_{|x|=1}2p(1-\alpha)\log(1-xz)d\mu(x)\right]\;\;z\in\mathbb{U}^{*}.
\end{equation}
The correspondence between $\Sigma_{\lambda\mu m p}(\alpha,1)$ and the probability measure $\mu(x)$ is one-to-one.
\end{theo}
\textit{Proof.} In view of the subordination condition (2.5), we have that $%
f\in\Sigma_{\lambda\mu m p}(\alpha,1)$ if and only if 
\begin{equation*}
-\frac{1}{p}\frac{z(D_{\lambda\mu p}^{m}f(z))^{\prime }}{D_{\lambda\mu
p}^{m}f(z)}\in\mathcal{P}(\alpha).
\end{equation*}
From Lemma 1.2, we have 
\begin{equation*}
-\frac{1}{p}\frac{z(D_{\lambda\mu p}^{m}f(z))^{\prime }}{D_{\lambda\mu
p}^{m}f(z)}=\int_{|x|=1}\frac{1-(2\alpha-1)xz}{1-xz}d\mu(x)
\end{equation*}
which is equivalent to 
\begin{equation*}
\frac{z(D_{\lambda\mu p}^{m}f(z))^{\prime }}{D_{\lambda\mu p}^{m}f(z)}%
=\int_{|x|=1}\frac{p(2\alpha-1)xz-p}{1-xz}d\mu(x).
\end{equation*}
Integrating this equality, we obtain 
\begin{equation*}
z^{p}D_{\lambda\mu p}^{m}f(z)=\exp\int_{|x|=1}2p(1-\alpha)\log(1-xz)d\mu(x)
\end{equation*}
or 
\begin{equation}  \label{2.7}
D_{\lambda\mu p}^{m}f(z)=z^{-p}\exp\int_{|x|=1}2p(1-\alpha)\log(1-xz)d\mu(x).
\end{equation}
Equality (2.6) now follows easily from (1.12), (1.13) and (2.7).

Using a result of Goluzin \cite{gol} (see also \cite{pom} p. 50), we obtain
the following result. 
\begin{theo}
Let $f\in\Sigma_{\lambda\mu m p}(\alpha,1)$. Then
$$z^{p}D_{\lambda\mu p}^{m}f(z)\prec(1-z)^{2p(1-\alpha)}\;\;(z\in\mathbb{U}).$$
\end{theo}
\textit{Proof.} Let $f\in\Sigma_{\lambda\mu m p}(\alpha,1)$. Then by (2.5)
we have 
\begin{equation*}
\frac{z(D_{\lambda\mu p}^{m}f(z))^{\prime }}{D_{\lambda\mu p}^{m}f(z)}\prec%
\frac{p(2\alpha-1)z-p}{1-z}.
\end{equation*}
Since the function $\displaystyle\frac{p(2\alpha-1)z-p}{1-z}$ is univalent
and convex in $\mathbb{U}$, in view of Goluzin's result, we obtain 
\begin{equation*}
\int_{0}^{z}\frac{(D_{\lambda\mu p}^{m}f(\zeta))^{\prime }}{D_{\lambda\mu
p}^{m}f(\zeta)}d\zeta\prec\int_{0}^{z}\frac{p(2\alpha-1)\zeta-p}{%
\zeta(1-\zeta)}d\zeta
\end{equation*}
or 
\begin{equation*}
\log(D_{\lambda\mu p}^{m}f(z))\prec\log\frac{(1-z)^{2p(1-\alpha)}}{z^{p}}.
\end{equation*}
Thus, there exists a function $w\in\Omega$ such that 
\begin{equation*}
\log(D_{\lambda\mu p}^{m}f(z))=\log\frac{(1-w(z))^{2p(1-\alpha)}}{w(z)^{p}}
\end{equation*}
which is equivalent to 
\begin{equation*}
z^{p}D_{\lambda\mu p}^{m}f(z)\prec(1-z)^{2p(1-\alpha)}.
\end{equation*}

Next we obtain a structural formula for the class $\Sigma_{\lambda\mu m
p}(\alpha,\beta)$. 
\begin{theo}
Let $f\in\Sigma_{\lambda\mu m p}(\alpha,\beta)$. Then 
$$f(z)=\left[z^{-p}+\sum_{k=1-p}^{\infty}\frac{z^{k}}{\Phi_{k}(\lambda,\mu,m,p)}\right]$$
\begin{equation}\label{2.8}
\ast\left[z^{-p}\exp\left(2p(1-\alpha)\beta\int_{0}^{z}\frac{w(\zeta)}{\zeta(1-\beta w(\zeta))}d\zeta\right)\right]\;(z\in\mathbb{U}^{*})
\end{equation}
where $w\in\Omega$.
\end{theo}
\textit{Proof.} Assume $f\in\Sigma_{\lambda\mu m p}(\alpha,\beta)$. From
(2.1) it follows 
\begin{equation}  \label{2.9}
\frac{z(D_{\lambda\mu p}^{m}f(z))^{\prime }}{D_{\lambda\mu p}^{m}f(z)}=\frac{%
p(2\alpha-1)\beta w(z)-p}{1-\beta w(z)}\;\;(z\in\mathbb{U}).
\end{equation}
In view of (2.9), we have 
\begin{equation*}
\frac{(D_{\lambda\mu p}^{m}f(z))^{\prime }}{D_{\lambda\mu p}^{m}f(z)}+\frac{p%
}{z}=\frac{p(2\alpha-1)\beta w(z)}{z(1-\beta w(z))}\;\;\;(z\in\mathbb{U}^{*})
\end{equation*}
which upon integration, yields 
\begin{equation}  \label{2.10}
\log(z^{p}D_{\lambda\mu p}^{m}f(z))=2p(1-\alpha)\beta\int_{0}^{z}\frac{%
w(\zeta)}{\zeta(1-\beta w(\zeta))}d\zeta.
\end{equation}
The assertion (2.8) of the theorem can be easily obtained from (1.12),
(1.13) and (2.10).

In the sequence a convolution property for the class $\Sigma_{\lambda\mu m
p}(\alpha,\beta)$ is derived. 
\begin{theo}
If $f\in\Sigma_{p}$ belongs to $\Sigma_{\lambda\mu m p}(\alpha,\beta)$, then
\begin{equation}\label{2.11}
D_{\lambda\mu p}^{m}f(z)\ast\left\{\frac{-pz^{-p}+(p+1)z^{-p+1}}{(1-z)^{2}}(1-\beta e^{i\theta})+\frac{z^{-p}}{1-z}[p-p(2\alpha-1)\beta e^{i\theta}]\right\}\neq0
\end{equation}
for $z\in\mathbb{U}^{*}$ and $\theta\in(0,2\pi)$.
\end{theo}
\textit{Proof.} Let $f\in\Sigma_{\lambda\mu m p}(\alpha,\beta)$. Then, from
(2.1) it follows 
\begin{equation}  \label{2.12}
-\frac{z(D_{\lambda\mu p}^{m}f(z))^{\prime }}{D_{\lambda\mu p}^{m}f(z)}\neq%
\frac{p-p(2\alpha-1)\beta e^{i\theta}}{1-\beta e^{i\theta}}\;\;(z\in\mathbb{U%
}\;,\;0<\theta<2\pi).
\end{equation}
It is easy to see that the condition (2.12) can be written as follows 
\begin{equation}  \label{2.13}
(1-\beta e^{i\theta})z(D_{\lambda\mu p}^{m}f(z))^{\prime
i\theta}]D_{\lambda\mu p}^{m}f(z)\neq0.
\end{equation}
Note that 
\begin{equation*}
D_{\lambda\mu p}^{m}f(z)=D_{\lambda\mu
p}^{m}f(z)\ast\left(z^{-p}+z^{-p+1}+\ldots+\frac{1}{z}+1+\frac{z}{1-z}\right)
\end{equation*}
\begin{equation}  \label{2.14}
=D_{\lambda\mu p}^{m}f(z)\ast\frac{z^{-p}}{1-z}
\end{equation}
and 
\begin{equation*}
z(D_{\lambda\mu p}^{m}f(z))^{\prime }=D_{\lambda\mu p}^{m}f(z)\ast\left[%
-pz^{-p}-(p-1)z^{-p+1}-\ldots-\frac{1}{z}+\frac{z}{(1-z)^{2}}\right]
\end{equation*}
\begin{equation}  \label{2.15}
=D_{\lambda\mu p}^{m}f(z)\ast\frac{-pz^{-p}+(p+1)z^{-p+1}}{(1-z)^{2}}.
\end{equation}
By virtue of (2.13), (2.14) and (2.15), the assertion (2.12) of the theorem
follows.

Coefficient estimates for functions in the class $\Sigma_{\lambda\mu m
p}(\alpha,\beta)$ are given in the next theorem. 
\begin{theo}
Let $f$ of the form (1.6) be in the class $\Sigma_{\lambda\mu m p}(\alpha,\beta)$. Then, for $n\geq 3-p$
\begin{equation}\label{2.16}
|a_{n}|\leq\frac{2p\beta(1-\alpha)}{(n+p)\Phi_{n}(\lambda,\mu,m,p)}
\end{equation}
where $\Phi_{n}(\lambda,\mu,m,p)$ is given by (1.11).
\end{theo}
\textit{Proof.} To prove the coefficient estimates (2.16) we use the method
of Clunie and Koegh \cite{clu}.

Let $f\in\Sigma_{\lambda\mu m p}(\alpha,\beta)$. From (1.14), we have 
\begin{equation*}
\frac{1}{p}\frac{z(D_{\lambda\mu p}^{m}f(z))^{\prime }}{D_{\lambda\mu
p}^{m}f(z)}+1=zw(z)\left[\frac{1}{p}\frac{z(D_{\lambda\mu
p}^{m}f(z))^{\prime }}{D_{\lambda\mu p}^{m}f(z)}+2\alpha-1\right]
\end{equation*}
where $w$ is analytic in $\mathbb{U}$ and $|w(z)|\leq\beta$ for $z\in\mathbb{%
U}$. Then 
\begin{equation}  \label{2.17}
z(D_{\lambda\mu p}^{m}f(z))^{\prime }+pD_{\lambda\mu p}^{m}f(z)=zw(z)\left[%
z(D_{\lambda\mu p}^{m}f(z))^{\prime }+p(2\alpha-1)D_{\lambda\mu p}^{m}f(z)%
\right].
\end{equation}
If $zw(z)=\displaystyle\sum_{k=1}^{\infty}w_{k}z^{k}$, making use of (1.10)
and (2.17), we obtain 
\begin{equation*}
\sum_{k=1-p}^{\infty}(k+p)\Phi_{k}(\lambda,\mu,m,p)a_{k}z^{k+p}
\end{equation*}
\begin{equation}  \label{2.18}
=\left[-2p(1-\alpha)\sum_{k=1-p}^{\infty}[k+p(2\alpha-1)]\Phi_{k}(\lambda,%
\mu,m,p)a_{k}z^{k+p}\right]\sum_{k=1}^{\infty}w_{k}z^{k}.
\end{equation}
Equating the coefficients in (2.18), we have 
\begin{equation*}
n\Phi_{n-p}(\lambda,\mu,m,p)a_{n-p}=-2p(1-\alpha)w_{n}\;,\,\text{for}\;n=1,2
\end{equation*}
and 
\begin{equation*}
n\Phi_{n-p}(\lambda,\mu,m,p)a_{n-p}
\end{equation*}
\begin{equation}  \label{2.19}
=-2p(1-\alpha)w_{n}+\sum_{k=1-p}^{n-1-p}[k+p(2\alpha-1)]\Phi_{k}(\lambda,%
\mu,m,p)a_{k}w_{n-p-k}
\end{equation}
for $n\geq3$.

From (2.19), we obtain 
\begin{equation*}
\left[-2p(1-\alpha)+\sum_{k=1-p}^{n-1-p}[k+p(2\alpha-1)]\Phi_{k}(\lambda,%
\mu,m,p)a_{k}z^{k+p}\right]\sum_{k=1}^{\infty}w_{k}z^{k}
\end{equation*}
\begin{equation}  \label{2.20}
=\sum_{k=1-p}^{n-p}(k+p)\Phi_{k}(\lambda,\mu,m,p)a_{k}z^{k+p}+%
\sum_{k=n+1-p}^{\infty}c_{k}z^{k+p}.
\end{equation}
It is known that, if $h(z)=\displaystyle\sum_{n=0}^{\infty}h_{n}z^{n}$ is
analytic in $\mathbb{U}$, then for $0<r<1$ 
\begin{equation}  \label{2.21}
\sum_{n=0}^{\infty}|h_{n}|^{2}r^{2n}=\frac{1}{2\pi}\int_{0}^{2\pi}|h(re^{i%
\theta})|^{2}d\theta.
\end{equation}
Since $\left|\displaystyle\sum_{k=1}^{\infty}w_{k}z^{k}\right|\leq\beta|z|<%
\beta$, making use of (2.20) and (2.21), we have 
\begin{equation*}
\sum_{k=1-p}^{n-p}(k+p)^{2}\Phi_{k}(\lambda,%
\mu,m,p)^{2}|a_{k}|^{2}r^{2(k+p)}+\sum_{k=n+1-p}^{%
\infty}|c_{k}|^{2}r^{2(k+p)}
\end{equation*}
\begin{equation*}
\leq\beta^{2}\left\{4p^{2}(1-\alpha)^{2}+\sum_{k=1-p}^{n-1-p}[k+p(2%
\alpha-1)]^{2}\Phi_{k}(\lambda,\mu,m,p)^{2}|a_{k}|^{2}r^{2(k+p)}\right\}.
\end{equation*}
Letting $r\rightarrow 1$, we obtain 
\begin{equation*}
\sum_{k=1-p}^{n-p}(k+p)^{2}\Phi_{k}(\lambda,\mu,m,p)^{2}|a_{k}|^{2}
\end{equation*}
\begin{equation*}
\leq
4p^{2}\beta^{2}(1-\alpha)^{2}+\sum_{k=1-p}^{n-1-p}\beta^{2}[k+p(2%
\alpha-1)]^{2}\Phi_{k}(\lambda,\mu,m,p)^{2}|a_{k}|^{2}.
\end{equation*}
The above inequality implies 
\begin{equation*}
n^{2}\Phi_{n-p}(\lambda,\mu,m,p)^{2}|a_{n-p}|^{2}
\end{equation*}
\begin{equation*}
\leq
4p^{2}\beta^{2}(1-\alpha)^{2}-\sum_{k=1-p}^{n-1-p}(k+p)^{2}\Phi_{k}(\lambda,%
\mu,m,p)^{2}|a_{k}|^{2}
\end{equation*}
\begin{equation*}
+\sum_{k=1-p}^{n-1-p}\beta^{2}[k+p(2\alpha-1)]^{2}\Phi_{k}(\lambda,%
\mu,m,p)^{2}|a_{k}|^{2}\leq 4p^{2}\beta^{2}(1-\alpha)^{2}.
\end{equation*}
Finally, replacing $n-p$ by $n$, we have 
\begin{equation*}
|a_{n}|\leq\frac{2p\beta(1-\alpha)}{(n+p)\Phi_{n}(\lambda,\mu,m,p)}.
\end{equation*}
Thus, the proof of our theorem is completed.

Theorem 2.6 enables us to obtain a distortion result for the class $%
\Sigma_{\lambda\mu m p}(\alpha,\beta)$. 
\begin{cor}
If $f\in\Sigma_{\lambda\mu m p}(\alpha,\beta)$ is given by (1.6), then for $0<|z|=r<1$
$$|f(z)|\geq\frac{1}{r^{p}}-2p\beta(1-\alpha)r^{1-p}\sum_{k=1-p}^{\infty}\frac{1}{(k+p)\Phi_{k}(\lambda,\mu,m,p)}$$
$$|f(z)|\leq\frac{1}{r^{p}}+2p\beta(1-\alpha)r^{1-p}\sum_{k=1-p}^{\infty}\frac{1}{(k+p)\Phi_{k}(\lambda,\mu,m,p)}$$
and
$$|f'(z)|\geq\frac{p}{r^{p+1}}-2p\beta(1-\alpha)r^{2-p}\sum_{k=1-p}^{\infty}\frac{k}{(k+p)\Phi_{k}(\lambda,\mu,m,p)}$$
$$|f'(z)|\leq\frac{p}{r^{p+1}}+2p\beta(1-\alpha)r^{2-p}\sum_{k=1-p}^{\infty}\frac{k}{(k+p)\Phi_{k}(\lambda,\mu,m,p)}.$$
\end{cor}

In the sequence we give a sufficient condition for a function to belong to
the class $\Sigma_{\lambda\mu m p}(\alpha,\beta)$. 
\begin{theo}
Let $f\in\Sigma_{p}$ be given by (1.6). If for $0\leq\alpha<1$ and $0<\beta\leq1$
\begin{equation}\label{2.22}
\sum_{k=1-p}^{\infty}[k(\beta+1)+p(1+\beta(2\alpha-1))]\Phi_{k}(\lambda,\mu,m,p)|a_{k}|\leq2p\beta(1-\alpha)
\end{equation}
then $f\in\Sigma_{\lambda\mu m p}(\alpha,\beta)$.
\end{theo}
\textit{Proof.} Assume that $f(z)=z^{-p}+\displaystyle\sum_{k=1-p}^{%
\infty}a_{k}z^{k}$. We have 
\begin{equation*}
M=\left|z(D_{\lambda\mu p}^{m}f(z))^{\prime }+pD_{\lambda\mu
p}^{m}f(z)\right|-\beta\left|z(D_{\lambda\mu p}^{m}f(z))^{\prime
}+p(2\alpha-1)D_{\lambda\mu p}^{m}f(z)\right|
\end{equation*}
\begin{equation*}
=\left|\sum_{k=1-p}^{\infty}(k+p)\Phi_{k}(\lambda,\mu,m,p)a_{k}z^{k}\right|
\end{equation*}
\begin{equation*}
-\beta\left|\frac{-2p(1-\alpha)}{z^{p}}+\sum_{k=1-p}^{\infty}[k+p(2%
\alpha-1)]\Phi_{k}(\lambda,\mu,m,p)a_{k}z^{k}\right|.
\end{equation*}
For $0<|z|=r<1$, we obtain 
\begin{equation*}
r^{p}M\leq\sum_{k=1-p}^{\infty}(k+p)\Phi_{k}(\lambda,\mu,m,p)|a_{k}|r^{p+k}
\end{equation*}
\begin{equation*}
-\beta\left[2p(1-\alpha)-\sum_{k=1-p}^{\infty}|k+p(2\alpha-1)|\Phi_{k}(%
\lambda,\mu,m,p)|a_{k}|r^{p+k}\right]
\end{equation*}
or 
\begin{equation*}
r^{p}M\leq\sum_{k=1-p}^{\infty}[k(\beta+1)+p(1+\beta(2\alpha-1))]\Phi_{k}(%
\lambda,\mu,m,p)|a_{k}|r^{p+k}-2p\beta(1-\alpha).
\end{equation*}
Since the above inequality holds for all $r\;(0<r<1)$, letting $r\rightarrow1
$, we have 
\begin{equation*}
M\leq\sum_{k=1-p}^{\infty}[k(\beta+1)+p(1+\beta(2\alpha-1))]\Phi_{k}(%
\lambda,\mu,m,p)|a_{k}|-2p\beta(1-\alpha).
\end{equation*}
Making use of (2.22), we obtain $M\leq0$, that is 
\begin{equation*}
\left|\frac{1}{p}\frac{z(D_{\lambda\mu p}^{m}f(z))^{\prime }}{D_{\lambda\mu
p}^{m}f(z)}+1\right|<\beta\left|\frac{1}{p}\frac{z(D_{\lambda\mu
p}^{m}f(z))^{\prime }}{D_{\lambda\mu p}^{m}f(z)}+2\alpha-1\right|.
\end{equation*}
Consequently, $f\in\Sigma_{\lambda\mu m p}(\alpha,\beta)$.

\section{Properties of the class $\Sigma_{\protect\lambda\protect\mu mp}^{+}(%
\protect\alpha,\protect\beta)$}

We begin this section by proving that the condition (2.22) is both necessary
and sufficient for a function to be in the class $\Sigma_{\lambda\mu
mp}^{+}(\alpha,\beta)$. 
\begin{theo}
Let $f\in\Sigma_{p}^{+}$. Then $f$ belongs to the class $\Sigma_{\lambda\mu mp}^{+}(\alpha,\beta)$ if and only if
$$\sum_{k=1-p}^{\infty}[k(\beta+1)+p(1+\beta(2\alpha-1))]\Phi_{k}(\lambda,\mu,m,p)a_{k}\leq2p\beta(1-\alpha).$$
\end{theo} 
\textit{Proof.} In view of Theorem 2.7, we have to prove "only if" part.
Assume that $f(z)=z^{-p}+\displaystyle\sum_{k=1-p}^{\infty}a_{k}z^{k}$, $%
a_{k}\geq0$ is in the class $\Sigma_{\lambda\mu mp}^{+}(\alpha,\beta)$. Then 
\begin{equation*}
\left|\frac{\displaystyle\frac{1}{p}\frac{z(D_{\lambda\mu
p}^{m}f(z))^{\prime }}{D_{\lambda\mu p}^{m}f(z)}+1}{\displaystyle\frac{1}{p}%
\frac{z(D_{\lambda\mu p}^{m}f(z))^{\prime }}{D_{\lambda\mu p}^{m}f(z)}%
+2\alpha-1}\right|=
\end{equation*}
\begin{equation*}
\left|\frac{\displaystyle\sum_{k=1-p}^{\infty}(k+p)\Phi_{k}(\lambda,%
\mu,m,p)a_{k}z^{k}}{\displaystyle\frac{2p(1-\alpha)}{z^{p}}-\displaystyle%
\sum_{k=1-p}^{\infty}[k+p(2\alpha-1)]\Phi_{k}(\lambda,\mu,m,p)a_{k}z^{k}}%
\right|<\beta
\end{equation*}
for all $z\in\mathbb{U}$. Since $\Re z\leq|z|$ for all $z$, it follows that 
\begin{equation}  \label{3.1}
\Re\left\{\frac{\displaystyle\sum_{k=1-p}^{\infty}(k+p)\Phi_{k}(\lambda,%
\mu,m,p)a_{k}z^{k}}{\displaystyle\frac{2p(1-\alpha)}{z^{p}}-\displaystyle%
\sum_{k=1-p}^{\infty}[k+p(2\alpha-1)]\Phi_{k}(\lambda,\mu,m,p)a_{k}z^{k}}%
\right\}<\beta.
\end{equation}
We choose the values of $z$ on the real axis such that $\frac{1}{p}\frac{%
z(D_{\lambda\mu p}^{m}f(z))^{\prime }}{D_{\lambda\mu p}^{m}f(z)}$ is real.
Upon clearing the denominator in (3.1) and letting $z\rightarrow1$ through
positive values, we obtain 
\begin{equation*}
\sum_{k=1-p}^{\infty}(k+p)\Phi_{k}(\lambda,\mu,m,p)a_{k}
\end{equation*}
\begin{equation*}
\leq2p\beta(1-\alpha)-\sum_{k=1-p}^{\infty}\beta[k+p(2\alpha-1)]%
\Phi_{k}(\lambda,\mu,m,p)a_{k}
\end{equation*}
or 
\begin{equation*}
\sum_{k=1-p}^{\infty}[k(\beta+1)+p(1+\beta(2\alpha-1))]\Phi_{k}(\lambda,%
\mu,m,p)a_{k}\leq2p\beta(1-\alpha).
\end{equation*}
Hence, the result follows. 
\begin{cor}
If $f\in\Sigma_{p}^{+}$, given by (1.7) is in the class $\Sigma_{\lambda\mu mp}^{+}(\alpha,\beta)$, then
\begin{equation}\label{3.2}
a_{n}\leq\frac{2p\beta(1-\alpha)}{[n(\beta+1)+p(1+\beta(2\alpha-1))]\Phi_{n}(\lambda,\mu,m,p)}\;,\;n\geq1-p
\end{equation}
with equality for the functions of the form
$$f_{n}(z)=\frac{1}{z^{p}}+\frac{2p\beta(1-\alpha)}{[n(\beta+1)+p(1+\beta(2\alpha-1))]\Phi_{n}(\lambda,\mu,m,p)}z^{n}.$$
\end{cor}

Coefficient estimates obtained in Corollary 3.1 enables us to give a
distortion result for the class $\Sigma_{\lambda\mu mp}^{+}(\alpha,\beta)$. 
\begin{theo}
If $f\in\Sigma_{\lambda\mu mp}^{+}(\alpha,\beta)$, then for $0<|z|=r<1$
$$|f(z)|\geq\frac{1}{r^{p}}-\frac{2p\beta(1-\alpha)}{[\beta(1-p+(2\alpha-1)p)+1]\Phi_{1-p}(\lambda,\mu,m,p)}r^{1-p}$$ and
$$|f(z)|\leq\frac{1}{r^{p}}+\frac{2p\beta(1-\alpha)}{[\beta(1-p+(2\alpha-1)p)+1]\Phi_{1-p}(\lambda,\mu,m,p)}r^{1-p}$$
where equality holds for the function
$$f_{p}(z)=\frac{1}{z^{p}}+\frac{2p\beta(1-\alpha)}{[\beta(1-p+(2\alpha-1)p)+1]\Phi_{1-p}(\lambda,\mu,m,p)}z^{1-p}$$ at $z=ir,r$.
\end{theo}
\textit{Proof.} Suppose $f\in\Sigma_{\lambda\mu mp}^{+}(\alpha,\beta)$.
Making use of inequality 
\begin{equation}  \label{3.3}
\sum_{k=1-p}^{\infty}a_{k}\leq\frac{2p\beta(1-\alpha)}{[\beta(1-p+(2%
\alpha-1)p)+1]\Phi_{1-p}(\lambda,\mu,m,p)}
\end{equation}
which follows easily from Theorem 3.1, the proof is trivial.

Now, we prove that the class $\Sigma_{\lambda\mu mp}^{+}(\alpha,\beta)$ is
closed under convolution. 
\begin{theo}
Let $h(z)=z^{-p}+\displaystyle\sum_{k=1-p}^{\infty}c_{k}z^{k}$ be analytic in $\mathbb{U}^{*}$ and $0\leq c_{k}\leq1$. If $f$ given by (1.7) is in the class $\Sigma_{\lambda\mu mp}^{+}(\alpha,\beta)$, then $f\ast h$ is also in the class $\Sigma_{\lambda\mu mp}^{+}(\alpha,\beta)$.
\end{theo}
\textit{Proof.} Since $f\in\Sigma_{\lambda\mu mp}^{+}(\alpha,\beta)$, then
by Theorem 3.1, we have 
\begin{equation*}
\sum_{k=1-p}^{\infty}[k(\beta+1)+p(1+\beta(2\alpha-1))]\Phi_{k}(\lambda,%
\mu,m,p)a_{k}\leq2p\beta(1-\alpha).
\end{equation*}
In view of the above inequality and the fact that 
\begin{equation*}
(f\ast h)(z)=z^{-p}+\sum_{k=1-p}^{\infty}a_{k}c_{k}z^{k}
\end{equation*}
we obtain 
\begin{equation*}
\sum_{k=1-p}^{\infty}[k(\beta+1)+p(1+\beta(2\alpha-1))]\Phi_{k}(\lambda,%
\mu,m,p)a_{k}c_{k}\leq
\end{equation*}
\begin{equation*}
\sum_{k=1-p}^{\infty}[k(\beta+1)+p(1+\beta(2\alpha-1))]\Phi_{k}(\lambda,%
\mu,m,p)a_{k}\leq2p\beta(1-\alpha).
\end{equation*}
Therefore, by Theorem 3.1, the result follows.

The next result involves an integral operator which was investigated in many
papers \cite{aou}, \cite{din}, \cite{ura}. 
\begin{theo}
If $f\in\Sigma_{\lambda\mu mp}^{+}(\alpha,\beta)$, then the integral operator
$$F_{c,p}(z)=\frac{c}{z^{p+c}}\int_{0}^{z}t^{c+p-1}f(t)dt\;,\;c>0$$
is also in the class $\Sigma_{\lambda\mu mp}^{+}(\alpha,\beta)$.
\end{theo}
\textit{Proof.} It is easy to ckeck that 
\begin{equation*}
F_{c,p}(z)=f(z)\ast\left(z^{-p}+\sum_{k=1-p}^{\infty}\frac{c}{c+p+k}%
z^{k}\right).
\end{equation*}
Since $0<\displaystyle\frac{c}{c+p+k}\leq1$, by Theorem 3.3, the proof is
trivial.

\section{Neighborhoods and partial sums}

Following earlier investigations on the familiar concept of neighborhoods of
analytic functions by Goodman \cite{goo}, Ruschweyh \cite{rus} and more
recently by Liu and Srivastava \cite{liu}, \cite{li}, Liu \cite{lli},
Altinta\c s et al. \cite{alt}, Orhan and Kamali \cite{ork}, Srivastava and
Orhan \cite{sro}, Orhan \cite{orh}, Deniz and Orhan \cite{den} and Aouf \cite%
{ao}, we define the $(n,\delta)-$ neighborhood of a function $f\in\Sigma_{p}$
of the form (1.6) as follows. 
\begin{defn}
For $\delta>0$ and a non-negative sequence $S=\left\{s_{k}\right\}^{\infty}_{k=1-p}$ where
\begin{equation}\label{4.1}
s_{k}:=\frac{[\beta(k+|2\alpha-1|p+k+p)]\Phi_{k}(\lambda,\mu, m,p)}{2p\beta(1-\alpha)}
\end{equation}
$$(k\geq1-p\;,\;p\in\mathbb{N}\;,\;0\leq\alpha<1\;,\;0<\beta\leq1)$$
the $(n,\delta)-$ neighborhood of a function $f\in\Sigma_{p}$ of the form (1.6) is defined by
\begin{equation}\label{4.2}
\mathcal{N}_{\delta}(f):=\left\{g\in\Sigma_{p}: g(z)=z^{-p}+\sum_{k=1-p}^{\infty}b_{k}z^{k}\;\textrm{and}\;\sum_{k=1-p}^{\infty}s_{k}|b_{k}-a_{k}|\leq\delta\right\}.
\end{equation}
\end{defn}

For $s_{k}=k$, Definition 1.4 corresponds to the $(n,\delta)-$ neighborhood
considered by Ruscheweyh \cite{rus}.

Making use of Definition 4.1, we prove the first result on $(n,\delta)-$
neighborhood of the class $\Sigma_{\lambda\mu mp}(\alpha,\beta)$. 
\begin{theo}
Let $f\in\Sigma_{\lambda\mu mp}(\alpha,\beta)$ be given by (1.6). If $f$ satisfies
\begin{equation}\label{4.3}
(f(z)+\epsilon z^{p})(1+\epsilon)^{-1}\in\Sigma_{\lambda\mu mp}(\alpha,\beta)\;\;(\epsilon\in\mathbb{C}\;,\;|\epsilon|<\delta\;,\;\delta>0),
\end{equation}
then
\begin{equation}\label{4.4}
\mathcal{N}_{\delta}(f)\subset\Sigma_{\lambda\mu mp}(\alpha,\beta).
\end{equation}
\end{theo}
\textit{Proof.} It is not difficult to see that a function $f$ belongs to $%
\Sigma_{\lambda\mu mp}(\alpha,\beta)$ if and only if 
\begin{equation}  \label{4.5}
\frac{z(D_{\lambda\mu p}^{m}f(z))^{\prime }+pD_{\lambda\mu p}^{m}f(z)}{\beta
z(D_{\lambda\mu p}^{m}f(z))^{\prime }+\beta(2\alpha-1)pD_{\lambda\mu
p}^{m}f(z)}\neq\sigma\;\;(z\in\mathbb{U},\sigma\in\mathbb{C},|\sigma|=1)
\end{equation}
which is equivalent to 
\begin{equation}  \label{4.6}
\frac{(f\ast h)(z)}{z^{-p}}\neq0\;\;(z\in\mathbb{U}),
\end{equation}
where for convenience, 
\begin{equation*}
h(z):=z^{-p}+\sum_{k=1-p}^{\infty}c_{k}z^{k}
\end{equation*}
\begin{equation}  \label{4.7}
=z^{-p}+\sum_{k=1-p}^{\infty}\frac{[\beta\sigma(k+(2\alpha-1)p)-(k+p)]%
\Phi_{k}(\lambda,\mu, m,p)}{2p\beta(1-\alpha)\sigma}z^{k}.
\end{equation}
From (4.7) it follows that 
\begin{equation*}
|c_{k}|=\left|\frac{[\beta\sigma(k+(2\alpha-1)p)-(k+p)]\Phi_{k}(\lambda,%
\mu,m,p)}{2p\beta(1-\alpha)\sigma}\right|
\end{equation*}
\begin{equation}  \label{4.8}
\leq\frac{[\beta\sigma(k+|2\alpha-1|p)+k+p]\Phi_{k}(\lambda,\mu,m,p)}{%
2p\beta(1-\alpha)\sigma}\;\;(k\geq1-p,p\in\mathbb{N}).
\end{equation}
Furthermore, under the hypotheses (4.3), using (4.6) we obtain the following
assertions: 
\begin{equation*}
\frac{((f(z)+\epsilon z^{p})(1+\epsilon)^{-1})\ast h(z)}{z^{-p}}%
\neq0\;\;(z\in\mathbb{U})
\end{equation*}
or 
\begin{equation*}
\frac{(f\ast h)(z)}{z^{-p}}\neq\epsilon\;\;(z\in\mathbb{U}),
\end{equation*}
which is equivalent to 
\begin{equation}  \label{4.9}
\left|\frac{(f\ast h)(z)}{z^{-p}}\right|\geq\delta\;\;(z\in\mathbb{U}%
,\delta>0).
\end{equation}
Now, if we let 
\begin{equation*}
g(z)=z^{-p}+\sum_{k=1-p}^{\infty}b_{k}z^{k}\in\mathcal{N}_{\delta}(f),
\end{equation*}
then we have 
\begin{equation*}
\left|\frac{(f(z)-g(z))\ast h(z)}{z^{-p}}\right|=\left|\sum_{k=1-p}^{%
\infty}(a_{k}-b_{k})z^{k+p}\right|
\end{equation*}
\begin{equation*}
\leq\sum_{k=1-p}^{\infty}\frac{[\beta(k+|2\alpha-1|p)+k+p]\Phi_{k}(\lambda,%
\mu,m,p)}{2p\beta(1-\alpha)}|a_{k}-b_{k}||z|^{k+p}<\delta
\end{equation*}
\begin{equation*}
(z\in\mathbb{U},\delta>0,k\geq1-p,p\in\mathbb{N}).
\end{equation*}
Thus, for any complex number $\sigma$ with $|\sigma|=1$, we have 
\begin{equation*}
\frac{(g\ast h)(z)}{z^{-p}}\neq0\;\;(z\in\mathbb{U})
\end{equation*}
which implies $g\in\Sigma_{\lambda\mu mp}(\alpha,\beta)$. The proof of the
theorem is completed.

In the sequence we give the definition of $(n,\delta)-$neighborhood of a
function $f\in\Sigma_{p}^{+}$ of the form (1.7). 
\begin{defn}
For $\delta>0$ and a non-negative sequence $S=\left\{s_{k}\right\}^{\infty}_{k=1-p}$ where
$$s_{k}:=\frac{[k(\beta+1)+p(1+\beta(2\alpha-1))]\Phi_{k}(\lambda,\mu, m,p)}{2p\beta(1-\alpha)}$$
$$(k\geq1-p\;,\;p\in\mathbb{N}\;,\;0\leq\alpha<1\;,\;0<\beta\leq1)$$
the $(n,\delta)-$ neighborhood of a function $f\in\Sigma_{p}^{+}$ of the form (1.7) is defined by
\begin{equation}\label{4.10}
\tilde{\mathcal{N}}_{\delta}(f):=\left\{g\in\Sigma_{p}^{+}: g(z)=z^{-p}+\sum_{k=1-p}^{\infty}b_{k}z^{k}\;\textrm{and}\;\sum_{k=1-p}^{\infty}s_{k}|b_{k}-a_{k}|\leq\delta\right\}.
\end{equation}
\end{defn}

We have the following result on $(n,\delta)-$neighborhood of the class $%
\Sigma_{\lambda\mu mp}^{+}(\alpha,\beta)$. 
\begin{theo}
If the function $f$ given by (1.7) is in the class $\Sigma_{\lambda\mu mp}^{+}(\alpha,\beta)$, then
\begin{equation}\label{4.11}
\tilde{\mathcal{N}}_{\delta}(f)\subset\Sigma_{\lambda\mu mp}^{+}(\alpha,\beta),
\end{equation}
where
$$\delta=\frac{2\lambda\mu+\lambda-\mu}{1+2\lambda\mu+\lambda-\mu}.$$
The result is the best possible in the sense that $\delta$ cannot be increased.
\end{theo}
\textit{Proof.} For a function $f\in\Sigma_{\lambda\mu mp}^{+}(\alpha,\beta)$
of the form (1.7), Theorem 3.1 immediately yields 
\begin{equation}  \label{4.12}
\sum_{k=1-p}^{\infty}\frac{[k(\beta+1)+p(1+\beta(2\alpha-1))]\Phi_{k}(%
\lambda,\mu, m,p)}{2p\beta(1-\alpha)}a_{k}\leq\frac{1}{\Phi_{1-p}(\lambda,%
\mu,1,p)}.
\end{equation}
Let 
\begin{equation*}
g(z)=z^{-p}+\sum_{k=1-p}^{\infty}b_{k}z^{k}\in\tilde{\mathcal{N}}_{\delta}(f)
\end{equation*}
with 
\begin{equation*}
\delta=\frac{2\lambda\mu+\lambda-\mu}{1+2\lambda\mu+\lambda-\mu}>0.
\end{equation*}
From the condition (4.10) we find that 
\begin{equation}  \label{4.13}
\sum_{k=1-p}^{\infty}s_{k}|b_{k}-a_{k}|\leq\delta.
\end{equation}
Using (4.12) and (4.13), we obtain 
\begin{equation*}
\sum_{k=1-p}^{\infty}s_{k}b_{k}\leq\sum_{k=1-p}^{\infty}s_{k}a_{k}+%
\sum_{k=1-p}^{\infty}s_{k}|b_{k}-a_{k}|
\end{equation*}
\begin{equation*}
\leq\frac{1}{\Phi_{1-p}(\lambda,\mu,1,p)}+\delta=1.
\end{equation*}
Thus, in view of Theorem 3.1, we get $g\in\Sigma_{\lambda\mu
mp}^{+}(\alpha,\beta)$.

To prove the sharpness of the assertion of the theorem, we consider the
functions $f\in\Sigma_{\lambda\mu mp}^{+}(\alpha,\beta)$ and $%
g\in\Sigma_{p}^{+}$ given by 
\begin{equation}  \label{4.14}
f(z)=z^{-p}+\frac{2p\beta(1-\alpha)}{[\beta(1-p+(2\alpha-1)p)+1]\Phi_{1-p}(%
\lambda,\mu, m,p)}z^{1-p}
\end{equation}
and 
\begin{equation*}
g(z)=z^{-p}+\left[\frac{2p\beta(1-\alpha)}{[\beta(1-p+(2\alpha-1)p)+1]%
\Phi_{1-p}(\lambda,\mu, m,p)} +\right.
\end{equation*}
\begin{equation}  \label{4.15}
+\left.\frac{2p\beta(1-\alpha)}{[\beta(1-p+(2\alpha-1)p)+1]\Phi_{1-p}(%
\lambda,\mu, m,p)}\delta^{*}\right]z^{1-p}
\end{equation}
where $\delta^{*}>\delta.$

Clearly, the function $g$ belongs to $\tilde{\mathcal{N}}_{\delta}(f)$ but
according to Theorem 3.1, $g\notin\Sigma_{\lambda\mu mp}^{+}(\alpha,\beta)$.
Consequently, the proof of our theorem is completed.

Next, we investigate the ratio of real parts of functions of the form (1.6)
and their sequences of partial sums defined by 
\begin{equation}  \label{4.16}
k_{m}(z)=\left\{%
\begin{array}{ll}
z^{-p}, & \mbox{  $m=1,2,\ldots,-p$} \\ 
z^{-p}+\displaystyle\sum_{k=1-p}^{m-1}a_{k}z^{k}, & \mbox{ 
$m=1-p,2-p,\ldots$}%
\end{array}
\right.
\end{equation}
We also determine sharp lower bounds for $\Re\left\{\displaystyle\frac{f(z)}{%
k_{m}(z)}\right\}$ and $\Re\left\{\displaystyle\frac{k_{m}(z)}{f(z)}\right\}$%
. 
\begin{theo}
Let $f\in\Sigma_{p}$ be given by (1.6) and let $k_{m}(z)$ be given by (4.16). Suppose that
\begin{equation}\label{4.17}
\sum_{k=1-p}^{\infty}\theta_{k}|a_{k}|\leq1
\end{equation}
where $$\theta_{k}=\frac{[k(\beta+1)+p(1+\beta(2\alpha-1))]\Phi_{k}(\lambda,\mu,m,p)}{2p\beta(1-\alpha)}.$$
Then, for $m\geq1-p$, we have
\begin{equation}\label{4.18}
\Re\left\{\frac{f(z)}{k_{m}(z)}\right\}>1-\frac{1}{\theta_{m}}
\end{equation}
and
\begin{equation}\label{4.19}
\Re\left\{\frac{k_{m}(z)}{f(z)}\right\}>\frac{\theta_{m}}{1+\theta_{m}}.
\end{equation}
The results are sharp for each $m\geq1-p$ with the extremal function given by
\begin{equation}\label{4.20}
f(z)=z^{-p}-\frac{1}{\theta_{m}}z^{m}.
\end{equation}
\end{theo}
\textit{Proof.} Under the hypotheses of the theorem, we can see from (4.17)
that 
\begin{equation*}
\theta_{k+1}>\theta_{k}>1\;\;(k\geq1-p).
\end{equation*}
Therefore, by using hypotheses (4.17) again, we have 
\begin{equation}  \label{4.21}
\sum_{k=1-p}^{m-1}|a_{k}|+\theta_{m}\sum_{k=m}^{\infty}|a_{k}|\leq%
\sum_{k=1-p}^{\infty}\theta_{k}|a_{k}|\leq1.
\end{equation}
Let

\begin{equation}  \label{4.22}
\omega(z)=\theta_{m}\left[\frac{f(z)}{k_{m}(z)}-\left(1-\frac{1}{\theta_{m}}%
\right)\right]=1+\frac{\theta_{m}\displaystyle\sum_{k=m}^{\infty}a_{k}z^{k+p}%
}{1+\displaystyle\sum_{k=1-p}^{m-1}a_{k}z^{k+p}}.
\end{equation}
Applying (4.21) and (4.22), we find 
\begin{equation*}
\left|\frac{\omega(z)-1}{\omega(z)+1}\right|=\left|\frac{\theta_{m}%
\displaystyle\sum_{k=m}^{\infty}a_{k}z^{k+p}}{2+2\displaystyle%
\sum_{k=1-p}^{m-1}a_{k}z^{k+p}+\theta_{m}\displaystyle\sum_{k=m}^{%
\infty}a_{k}z^{k+p}}\right|
\end{equation*}
\begin{equation}  \label{4.23}
\leq\frac{\theta_{m}\displaystyle\sum_{k=m}^{\infty}|a_{k}|}{2-2\displaystyle%
\sum_{k=1-p}^{m-1}|a_{k}|-\theta_{m}\displaystyle\sum_{k=m}^{\infty}|a_{k}|}%
\leq1\;\;(z\in\mathbb{U})
\end{equation}
which shows that $\Re\omega(z)>0\;(z\in\mathbb{U})$. From (4.22), we
immediately obtain (4.18).

To prove that the function $f$ defined by (4.20) gives sharp result, we can
see that for $z\rightarrow1^{-}$ 
\begin{equation*}
\frac{f(z)}{k_{m}(z)}=1-\frac{1}{\theta_{m}}z^{m}\rightarrow1-\frac{1}{%
\theta_{m}}
\end{equation*}
which shows that the bound in (4.18) is the best possible.

Similarly, if we let 
\begin{equation}  \label{4.24}
\phi(z)=(1+\theta_{m})\left[\frac{k_{m}(z)}{f(z)}-\frac{\theta_{m}}{%
1+\theta_{m}}\right]=1-\frac{(1+\theta_{m})\displaystyle\sum_{k=m}^{%
\infty}a_{k}z^{k+p}}{1+\displaystyle\sum_{k=1-p}^{\infty}a_{k}z^{k+p}}
\end{equation}
and making use of (4.21), we find that 
\begin{equation*}
\left|\frac{\phi(z)-1}{\phi(z)+1}\right|=\left|\frac{-(1+\theta_{m})%
\displaystyle\sum_{k=m}^{\infty}a_{k}z^{k+p}}{2+2\displaystyle%
\sum_{k=1-p}^{\infty}a_{k}z^{k+p}-(1+\theta_{m})\displaystyle%
\sum_{k=m}^{\infty}a_{k}z^{k+p}}\right|
\end{equation*}
\begin{equation}  \label{4.25}
\leq\frac{(1+\theta_{m})\displaystyle\sum_{k=m}^{\infty}|a_{k}|}{2-2%
\displaystyle\sum_{k=1-p}^{\infty}|a_{k}|+(1+\theta_{m})\displaystyle%
\sum_{k=m}^{\infty}|a_{k}|}\leq1\;\;(z\in\mathbb{U})
\end{equation}
which leads immediately to the assertion (4.19) of the theorem.

The bound in (4.19) is sharp for each $m\geq1-p$, with the extremal function 
$f$ given by (4.20). The proof of the theorem is now completed.\newline
\newline
\textbf{Acknowledgements}\newline
The present investigation was supported by Ataturk University Rectorship
under BAP Project (The Scientific and Research Project of Ataturk
University) Project no: 2010/28.

\pagebreak {\footnotesize \textit{Department of Mathematics, Faculty of
Science and Arts, \newline
Ataturk University, 25240, Erzurum Turkey\newline
E-mail}: horhan@atauni.edu.tr\newline
\newline
\textit{Faculty of Mathematics and Computer Science, "Transilvania"
University of Bra\c sov, 50091, Iuliu Maniu, 50, Bra\c sov, Romania\newline
E-mail}: dorinaraducanu@yahoo.com\newline
\newline
\textit{Department of Mathematics, Faculty of Science and Arts, \newline
Ataturk University, 25240, Erzurum Turkey\newline
E-mail}: edeniz36@yahoo.com }

\end{document}